\documentclass[10pt]{article}

\usepackage[utf8]{inputenc}
\usepackage{amssymb}
\usepackage{amsmath}
\usepackage{graphicx}
\newtheorem{lemma}{Lemma}
\newtheorem{theorem}{Theorem}
\newtheorem{remark}{Remark}

\usepackage{tikz}
\usetikzlibrary{shapes,arrows}
\tikzset{cblue/.style={circle, draw, thin,fill=green!20, scale=0.7}}
\tikzset{cred/.style={circle, draw, thin, fill=red!20, scale=0.7}}
\tikzset{cgreen/.style={circle, draw, thin, fill=cyan!20, scale=0.7}}
\tikzset{cblack/.style={circle, draw, thin, fill=black, scale=0.3}}
\date{}
\author{Lerner~E.Yu.}
\title{Instances of small size with no weakly stable matching for three-sided problem with complete cyclic preferences}
\begin{document}
\maketitle

\begin{abstract}
Given $n$ men, $n$ women, and $n$ dogs, we assume that each man has a complete preference list of women, while each woman does a complete preference list of dogs, and each dog does a complete preference list of men.  
We study the so-called 3D-CYC problem, i.e., a three-dimensional problem with cyclic preferences. 
We understand a matching as a collection of $n$ nonintersecting triples, each of which contains a man, a woman, and a dog. A matching is said to be nonstable, if one can find a man, a woman, and a dog, which belong to different triples and prefer each other to their current partners in the corresponding triples. Otherwise, the matching is said to be stable. 
According to the conjecture proposed by Eriksson, S\"ostrand, and Strimling (2006), the problem of finding a stable matching (the problem 3DSM-CYC) always has a solution. 
However, Lam and Paxton have proposed an algorithm for constructing preference lists in 3DSM-CYC of size $n=90$, 
which has allowed them to disprove the mentioned conjecture. The question on the existence of counterexamples of a lesser size remained open. 
The main value of this paper consists in reducing the size of the counterexample to~$n=20$. 
At the end part of the paper, we discuss a new variant of 3DSM, whose solution always exists. 
\end{abstract}

\section{Introduction}
This paper has issued from our attempts to reformulate in terms of the graph theory the C.-K.~Lam and C.G.~Plaxton results~\cite{Lam}, which disprove the conjecture proposed by K.~Eriksson, J.~S\"ostrand, and P.~Strimling~\cite{Eriksson}. Such a formulation would make the mentioned results more demonstrative. Our paper is the second one (after~\cite{old}) in the series of works devoted to this topic. In the framework of the considered problem, we can treat the paper~\cite{old} as a preparatory one; it solves a separate problem stated by P.~Bir\'o and E.~McDermid~\cite{Biro}. Below we give all the necessary definitions and introductory remarks; we also prove the result from~\cite{old} that we use here. 
We tried to make this paper understandable even for an unprepared reader.

Assume that there are $n$ men and $n$ women, and each one among them has a (incomplete) preference list of representatives of the opposite sex. We understand a stable matching as a partition of the set of men and women into heterosexual families, where there exists no tuple consisting of a man and a woman, who would prefer each other to their partners (if any). Initially, complete preference lists were studied by D.~Gale and L.S.~Shaple; in~\cite{Gale}, they prove that a stable matching necessarily exists and propose an $O(n^2)$-hard algorithm for forming it. Note that the mentioned algorithm is also applicable for the case of incomplete preference lists, but in such a case, some men and women, possibly, remain single (single people of opposite sex do not want to marry each other and cannot take partners from married persons). 
In particular, this algorithm is used in many countries in enrolling applicants to universities. In 2012, one mentioned the practice of using this algorithm, when presenting the Nobel Prize in Economics to Lloyd Shapley and Alvin Roth.
A certain modification of the Gale--Shaple algorithm for incomlete lists case (see, e.g.,~\cite{knuth}) allows one to find, if possible, a matching without single men and women or to prove its absence, otherwise. 

In~\cite{knuth}, D.E.~Knuth states the question whether it is possible to generalize the theory of stable matchings to the case of three genders. The most interesting variant in the $k$-gender case occurs when preferences are cyclic, i.e., representatives of the 1st gender rank representatives of the 2nd one, the latter do representatives of the 3rd gender, etc., and each representative of the $k$th gender has a preference list of representatives of the 1st one.

A tuple containing exactly one representative of each gender is called a family, and a set of nonintersecting families is called a matching. 
A matching is said to be nonstable, if one can find one representative of each gender would become ``happier'' when forming a new family. Otherwise, the matching is said to be a weakly stable.
In what follows, for brevity, we use the term ``a stable matching'' instead of the term ``a weakly stable matching''.

Let the number of representatives of each gender equal~$n$. The authors of~\cite{gurvich} prove that in the case of complete preference lists, a stable matching always exists, provided that $n\leqslant k$ (here $k$ is the number of genders). In~\cite{Eriksson}, K.~Eriksson et al. generalize this result for the case when $k=3$ and $n=k+1=4$. Ibid, they state the conjecture that the problem of finding a stable matching in the 3-gender case with complete preference lists 
(three-dimensional stable matching with cyclic preferences problem --- 3DSM-CYC or just 3DSM) has a solution with any~$n$. Using the statement of the satisfiability problem and performing an extensive computer-assisted search, the authors of~\cite{new} prove the validity of the conjecture stated by K.~Eriksson et al. for~$n=5$. In~\cite{Pittel2}, B.~Pittel proves that with random preference lists the mean value of stable matchings in 3DSM grows as $\Omega(n^2 \ln^2(n))$.

The paper~\cite{Biro} by P.~Bir\'o and E.~McDermid is devoted to studying 3DSMI (the 3-gender stable matching problem with incomplete preference lists). Note that in this problem, similarly to the two-gender case, some agents, possibly, remain single. But in contrast to the two-gender case, according to results obtained in~\cite{Biro}, 3DSMI is not necessarily solvable; P.~Bir\'o and E.~McDermid give an explicit example of 3DSMI of size $n=6$, where no stable matching exists. Moreover, they prove that the problem of establishing the solvability of 3DSMI is NP-complete.

The conjecture stated by K.~Eriksson et al. has been recently disproved by C.-K.~Lam and C.G.~Plaxton~\cite{Lam}. They associate 3DSMI with a certain 3DSM problem, where $n$ is 15 times greater than the initial size; this problem is solvable if and only if so is the initial 3DSMI. Therefore, the problem of establishing the solvability of 3DSM is NP-complete. 
Note that C.-K.~Lam and C.G.~Plaxton initially considered the case of arbitrary $k>2$ (the factor 15 is a specification of the factor $k ((k-1)^2+1))$ for $k=3$).

The example proposed in the paper~\cite{Biro} allows one to construct an instance of 3DSM with no stable matching for $n=90=6\times15$. 
The question on the existence of counterexamples of a lesser size remained open. 

An evident way to reduce the size of counterexamples of 3DSM is to solve the P.~Bir\'o and E.~McDermid problem that implies the search of instances of 3DSMI with no stable matching for $n<6$. We solve the mentioned problem in the paper~\cite{old}. We prove the absence of such instances for $n<3$ and construct several counterexamples for 3DSMI with $n=3$. Therefore, the result obtained in~\cite{Lam} allows one to construct an instance of 3DSM with no stable matching for $n=45$.

Another approach implies the reduction of the multiplier (its current value equals 15) when constructing an unsolvable instance of 3DSM from an unsolvable 3DSMI problem. In this paper, we prove that one can associate each instance of 3DSMI with no stable matching of size $n$ with an instance of 3DSM with the same property of size $8n$. This allows one to reduce the size of the counterexample of 3DSM to $3\times 8 = 24$.

Finally, one can use constructions, which associate unsolvable 3DSMI problems with unsolvable 3DSM problems only in certain particular cases. 
Developing constructions mentioned in the previous paragraph and making use of specific features of a certain concrete instance of 3DSMI with no stable matching, we have succeeded in proposing a counterexample for 3DSM of size $n=20$.

The rest part of the paper has the following structure. In Sect.~2, we present formal statements of the 3DSM-CYC and the 3DSMI-CYC in terms of the graph theory. In Sect.~3, we prove that one can associate each instance of the 3DSMI of size $n$ with no stable matching with an instance of the 3DSM of size $8n$ with the same property. 
In Sect.~4, we improve the general scheme for further decreasing the size of counterexamples. 
In Sect.~5, we propose a counterexample of the 3DSM of size~20. In Sect.~6, we summarize the obtained results and discuss some potential future work. In particular, we propose a new variant of 3DSM, whose solution always exists and one can find it within quadratic time. In the Appendix, we give a short proof of Theorem~2 from~\cite{old}; we use it in this paper.

\section{The statement of 3DSM (3DSMI) in terms of the graph theory}
Let $G$ be some directed graph. Denote the set of its edges by $E$ (or $E(G)$); assume that no edge is multiple. Assume that the vertex set~$V$ of the graph~$G$ is divided into three subsets, namely, the set of men~$M$, women~$F$, and dogs~$D$. 
Assume that edges $(v,v')$, $v,v'\in V$, of this graph are such that either $v\in M, v'\in F$, or $v\in F, v'\in D$, or $v\in D, v'\in M$.
Assume that $|M|=|F|=|D|$ (otherwise we supplement the corresponding subgraph with vertices that are not connected with the rest part of the graph). The number $n=|M|=|F|=|D|$ is called the problem size. Evidently, the length of all cycles in the graph~$G$ is a multiple of~$3$. Note also that this condition ensures the possibility of dividing the vertex set of any digraph~$G$ into three subsets $M$, $F$, and $D$ so that all its edges are directed as is indicated above.

Each edge $(v,v')$, $v,v'\in V$, corresponds to some positive integer $r(v,v')$; it is called the rank of this edge. For fixed $v\in V$, all possible ranks $r(v,v_1),\ldots,r(v,v_k)$ coincide with $\{1,\ldots,k\}$, where $k$ is the outgoing vertex degree~$v$ (if $r(v, v') = 1$, then $v'$ is the best preference for $v$, and so on).

We understand a \textit{three-sided matching} as a subgraph~$H$ of the graph $G$, $V(H)=V(G)=V$, where each vertex $v\in V$ has at most one outgoing edge and the following condition is fulfilled: if a vertex $v$ has an outgoing edge, then this edge belongs to a cycle of length~3 in the graph~$H$. Cycles of length~3 in the graph $H$ are called families. Evidently, each family, accurate to a cyclic shift, takes the form $(m,f,d)$, where $m\in M$, $f\in F$, and $d\in D$. Note that in what follows, for convenience of denotations of families, we do not fix the order of genders in a family, i.e., we treat denotations of families as triples derived from an initial one by a cyclic shift as equivalent.

In what follows, we sometimes use the notion of a family in a wider sense, namely, as any cycle of length~3 in the graph~$G$. However, if some three-sided matching~$H$ is fixed, then we describe other cycles of length~3 explicitly, applying the term ``a family'' only to cycles that enter in a three-sided matching.

A \textit{matching} $\mu$ is a collection of all families of a three-sided matching $H$. For a vertex $v$, $v\in V$, in the matching $\mu$, the rank $R_\mu(v)$ is defined as the rank of the edge that goes out of this vertex in the subgraph $H$. If some vertex $v$ in the subgraph $H$ has no outgoing edge, then $R_\mu(v)$ is set to $+\infty$.

A triple $(v,v',v'')$ is said to be \textit{blocking} for some matching $\mu$, if it represents a cycle in the graph~$G$, and
\begin{equation}
\label{oldc}
r(v,v')<R_\mu(v),\quad r(v',v'')<R_\mu(v'),\quad r(v'',v)<R_\mu(v'').
\end{equation}
A matching $\mu$ is said to be \textit{stable}, if no blocking triple exists for it.

Recall that \textit{3DSMI} consists in finding a stable matching for a given graph~$G$. 
As is well known, such a matching does not necessarily exists. Moreover, the problem of establishing its existence for a given graph~$G$ is NP-complete. As was mentioned in the Introduction, the proof of this fact belongs to P.~Bir\'o and E.~McDermid. They have constructed an explicit example of the graph~$G$ of size~6, for which no stable matching exists. 

Evidently, 3DSM represents a particular case of the 3DSMI, where the outgoing (and incoming) degree of each vertex of the corresponding graph equals the problem size~$n$.

Let us construct on the set of agents (graph vertices) the same map as that used in~\cite{manlove} and in other papers. Denote it by the symbol~$\mu$; the same symbol stands for the matching, from which we construct this map. If an agent $x$ remains single, then we put $\mu(x)=x$. Otherwise $\mu(x)=y$, where $y$ is the vertex, which represents the endpoint of an edge in the subraph~$H$ that generates the matching~$\mu$. Informally speaking, in this case, $\mu(x)$ is the partner of $x$ in the family, whom~$x$ is ``partial'' to. Evidently, in 3DSM, as distinct from 3DSMI, the equality $\mu(x)=x$ is impossible for a stable matching.

In what follows, we use symbols $G$ and $H$ for graphs of 3DSM and 3DSMI, as well as their subgraphs; we also use, when appropriate, the symbol $H$ with various subscripts. At the same time, we distinguish between the subgraph $H$ and subgraphs denoted by the symbol $H$ with various subscripts, when considering them concurrently. In such cases, the symbol $H$ denotes the ``central'', in a sense, subgraph, while the symbol $H$ with some subscript does a ``peripheral'' subgraph characterized by its subscript. We use the same subscript in denotations for vertices of the latter subgraph. 

\section{The correspondence between unsolvable 3DSMI and 3DSM}
In this section, we prove that each 3DSMI of size~$n$ with no stable matching corresponds to 3DSM of size~$8n$ with the same property. As was mentioned in the Introduction, an analogous result with the multiplier~15 belongs to C.-K.~Lam and C.G.~Plaxton~\cite{Lam}. As a corollary, making use of results obtained by us in the paper~\cite{old} (see Lemma~\ref{ap} in the Appendix), we obtain concrete instances of 3DSM of size~24 with no stable matching. 

Let us consider subgraphs of the weighted graph of 3DSM, which include some vertices and edges of the initial graph. Let us first prove the Key Lemma~\ref{key1}.
(We divide the proof into several separate parts). 
We ``attach'' one copy of the subgraph considered in this lemma to each vertex of the graph~$H$ that defines 3DSMI with no stable matching. 

\begin{figure}[h]
\begin{center}
        \begin{tikzpicture}[->,>=stealth',shorten >=1pt,auto,node distance=3cm,thick,main node/.style={rectangle,fill=blue!20,draw,font=\sffamily\Large\bfseries}]
        \node[cred] (c) at (0:0) {$\mathbf c$};
        \node[cred] (d) at (0:4) {$\mathbf d$};
        \node[cgreen] (e) at (90:1) {$e$};
        \node[cblue] (a) at ( 330:1.6) {$a$};
        \node[cblue] (t) at ( 29:3.3) {$t$};
        \node[cgreen] (x) at ( 292.5:1.84776) {$x$};
        \node[cblue] (b) at ( 225:1.2) {$b$};
        \node[cgreen] (s) at (15:3) {$s$};

        \path[every node/.style={font=\sffamily\small}]
        (e) edge  []  node [] {} (c)
        (c) edge  []  node [] {} (a)
        (d) edge  []  node [] {} (a)
        (a) edge  []  node [] {} (x)
        (b) edge  []  node [] {} (x)
        (t) edge  [dashed]  node [] {} (s)
        (t) edge  []  node [] {} (e)
        (s) edge  []  node [] {} (d)
        (e) edge  [dashed]  node [] {} (d)
        (c) edge  [dashed]  node [right] {} (b)
        (d) edge  [dashed]  node [right] {} (b)
        (x) edge  [dotted]  node [right] {$r'_x$} (c)
        (x) edge  [dotted]  node [right] {$r_x'+1$} (d)    
        (b) edge  [dotted]  node [right] {3} (s)    
        (d) edge  [dotted]  node [right] {3} (t)    
        (a) edge  [dashed]  node [] {} (e)
        (b) edge  [dashed]  node [] {} (e)
        ;
        \end{tikzpicture}
\caption{The part of the preference graph considered in Lemma~\ref{key1}. 
Vertices of various colors correspond to various genders. Bold lines represent edges of rank~1, while dashed ones do those of rank~2. Ranks of edges represented by dotted lines equal $r'_x$, $r'_x+1$, or~3 (ranks are indicated near edges).}
\label{keyGraph}
\end{center}
\end{figure}
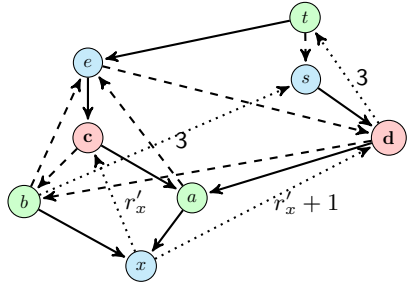

\begin{lemma}[The Key Lemma for Theorem~\ref{unSolvabable}]
\label{key1}
Let some subgraph of the graph of 3DSM take the form shown in~Fig.~\ref{keyGraph}, in particular,
\begin{equation}
\label{ranksRestrict}
r(b,s)=r(d,t)=3,\quad r(x,c)=r'_x,\ r(x,d)=r'_x+1,
\end{equation}
where $r'_x\in\{1,\ldots,n-1\}$. Let $\mu$ be a stable matching in this problem, while $y=\mu(x)$. Then only one of the following two alternatives is possible:\\ 
A) $y\in\{c,d\}$ and $\mu(y)\in\{a,b,t\}$; \\
B) $r(x,y)<r'_x$.
\end{lemma}

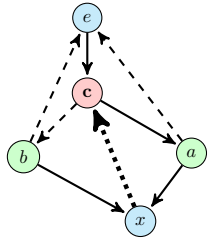
\begin{figure}[h]
\begin{center}
        \begin{tikzpicture}[->,>=stealth',shorten >=1pt,auto,node distance=3cm,thick,main node/.style={rectangle,fill=blue!20,draw,font=\sffamily\Large\bfseries}]
        \node[cred] (c) at (0:0) {$\mathbf c$};
        \node[cgreen] (e) at (90:1) {$e$};
        \node[cblue] (a) at ( 330:1.6) {$a$};
        \node[cgreen] (x) at ( 292.5:1.84776) {$x$};
        \node[cblue] (b) at ( 225:1.2) {$b$};

        \path[every node/.style={font=\sffamily\small}]
        (e) edge  []  node [] {} (c)
        (c) edge  []  node [] {} (a)
        (a) edge  []  node [] {} (x)
        (b) edge  []  node [] {} (x)
        (c) edge  [dashed]  node [right] {} (b)
        (x) edge  [line width=0.6mm, dotted]  node [left] {} (c)
        (a) edge  [dashed]  node [] {} (e)
        (b) edge  [dashed]  node [] {} (e)
        ;
        \end{tikzpicture}
\caption{The part of the preference graph considered in Lemma~\ref{key11}.
}
\label{keyGraph11}
\end{center}
\end{figure}

\begin{lemma}\label{key11}
Assume that some subgraph of the graph of 3DSM takes the form shown in Fig.~\ref{keyGraph11}.
Let $\mu$ be a stable matching in this problem, $\mu(x)=c$. Then $\mu(c)\in\{a,b\}$.
\end{lemma}

\noindent
\textbf{Proof of Lemma~\ref{key11}:}
Note that by assumption, $R_\mu(e)>1$. Assume that $\mu(c)\not\in\{a,b\}$. Then $R_\mu(c)>2$, $R_\mu(a)>1$, and $R_\mu(b)>1$. 
Consider the triple $(c,a,e)$. We get inequalities $r(c,a)<R_\mu(c)$, $r(e,c)<R_\mu(e)$, and $r(a,e)\leqslant R_\mu(a)$. Therefore, the triple $(c,a,e)$ is not blocking, only if $\mu(a)=e$. But then $R_\mu(b)>2$ and, consequently, the triple $(c,b,e)$ is blocking. 
\qquad  $\square$

In what follows, we repeatedly apply the technique that is used in the proof of this lemma. Namely, when considering the potentially blocking triple $(c,a,e)$, with the help of the mentioned technique we conclude that $\mu(a)=e$.

Note also that by the definition of a 3D-matching, for any distinct vertices $h$ and~$g$, equalities $\mu(h)=g$ and $\mu(\mu(g))=h$ are equivalent.
In other words, if $\mu(u_2)=u_1$, then the triple $(u_1, v,u_2)$ does not form families that enter in the matching, only if $\mu(u_1)\neq v$.

\begin{lemma}
\label{key111}
Assume that the graph of 3DSM has the subgraph shown in Fig.~\ref{keyGraph}.
Let $\mu$ be a stable matching in this problem, while $\mu(x)=d$. Then $\mu(d)\in\{a,b,t\}$.
\end{lemma}

\noindent
\textbf{Proof of Lemma~\ref{key111}:}
First of all, note that by assumption, $R_\mu(x)=r'_x+1$, $\mu(e) \neq d$, and $\mu(s)\neq d$. 

Let us prove the desired assertion ab contrario. Assume that $\mu(d)\not\in\{a,b,t\}$. 
Then $R_\mu(d)>3$. Moreover, $\mu(a) \neq x$, otherwise $\mu(d) = \mu(\mu(x))=a$. Therefore, $R_\mu(a)>1$. Analogously, $R_\mu(b)>1$.

Consequently, $\mu(c) =a$, otherwise the triple $(c,a,x)$ is blocking. 

Let us prove that $\mu(a)=e$, i.e., the matching $\mu$ contains the family $(c,a,e)$. Really, otherwise $R_\mu(a)>2$, $R_\mu(e)>1$, and the triple $(d,a,e)$ is blocking. 

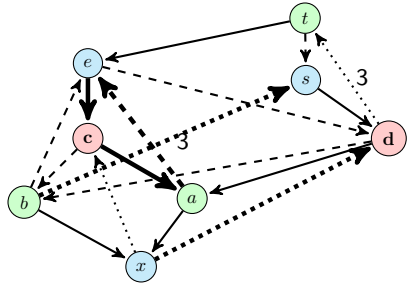
\begin{figure}[h]
\begin{center}
        \begin{tikzpicture}[->,>=stealth',shorten >=1pt,auto,node distance=3cm,thick,main node/.style={rectangle,fill=blue!20,draw,font=\sffamily\Large\bfseries}]
        \node[cred] (c) at (0:0) {$\mathbf c$};
        \node[cred] (d) at (0:4) {$\mathbf d$};
        \node[cgreen] (e) at (90:1) {$e$};
        \node[cblue] (a) at ( 330:1.6) {$a$};
        \node[cblue] (t) at ( 29:3.3) {$t$};
        \node[cgreen] (x) at ( 292.5:1.84776) {$x$};
        \node[cblue] (b) at ( 225:1.2) {$b$};
        \node[cgreen] (s) at (15:3) {$s$};

        \path[every node/.style={font=\sffamily\small}]
        (e) edge  [line width=0.6mm]  node [] {} (c)
        (c) edge  [line width=0.6mm]  node [] {} (a)
        (d) edge  []  node [] {} (a)
        (a) edge  []  node [] {} (x)
        (b) edge  []  node [] {} (x)
        (t) edge  [dashed]  node [] {} (s)
        (t) edge  []  node [] {} (e)
        (s) edge  []  node [] {} (d)
        (e) edge  [dashed]  node [] {} (d)
        (c) edge  [dashed]  node [right] {} (b)
        (d) edge  [dashed]  node [right] {} (b)
        (x) edge  [dotted]  node [] {} (c)
        (x) edge  [line width=0.6mm,dotted]  node [] {} (d)    
        (b) edge  [line width=0.6mm,dotted]  node [right] {3} (s)    
        (d) edge  [dotted]  node [right] {3} (t)    
        (a) edge  [line width=0.6mm,dashed]  node [] {} (e)
        (b) edge  [dashed]  node [] {} (e)
        ;
        \end{tikzpicture}
\caption{The part of the preference graph considered in Lemma~\ref{key111}. Bold edges correspond to pairs that enter in families of the matching $\mu$, provided that $\mu(d)\not\in\{a,b,t\}$.}
\label{keyGraphN}
\end{center}
\end{figure}

Analogously, $\mu(b) =s$ (otherwise the triple $(d,b,s)$ is blocking) (see Fig.~\ref{keyGraphN}).

Since $\mu(t)\not\in\{e,s\}$, we get the inequality $R_\mu(t)>2$. But then the triple $(d,t,s)$ is blocking.
\qquad  $\square$

\begin{figure}[h]
\begin{center}
        \begin{tikzpicture}[->,>=stealth',shorten >=1pt,auto,node distance=3cm,thick,main node/.style={rectangle,fill=blue!20,draw,font=\sffamily\Large\bfseries}]
        \node[cred] (c) at (0:0) {$\mathbf c$};      
        \node[cred] (d) at (0:4) {$\mathbf d$};
        \node[cblue] (a) at ( 330:1.6) {$a$};
        \node[cgreen] (x) at ( 292.5:1.84776) {$x$};
        \node[cblue] (b) at ( 225:1.2) {$b$};   
        
        \path[every node/.style={font=\sffamily\small}]
        (c) edge  []  node [] {} (a)
        (d) edge  []  node [] {} (a)
        (a) edge  [line width=0.6mm]  node [] {} (x)
        (b) edge  []  node [] {} (x)
        (c) edge  [line width=0.6mm,dashed]  node [right] {} (b)
        (d) edge  [dashed]  node [right] {} (b)
        (x) edge  [dotted]  node [right] {$r'_x$} (c)
        (x) edge  [dotted]  node [right] {$r_x'+1$} (d)    
        ;
        \end{tikzpicture}
\caption{The part of the preference graph considered in Lemma~\ref{key13}. 
Solid lines represent edges of rank~1, while dashed ones do those of rank~2. 
Ranks of edges represented by dotted lines equal $r'_x$ and $r'_x+1$ (ranks are indicated near edges).
Bold edges correspond to pairs that enter in families of the matching $\mu$ in the case, when $R_\mu(x)\geqslant r'(x)$, $\mu(c)\neq a$.}
\label{keyGraphOld}
\end{center}
\end{figure}
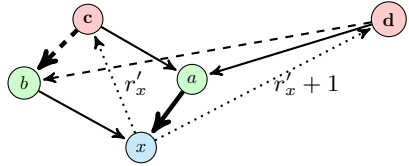

\begin{lemma}
\label{key13}
Assume that some subgraph of the graph of 3DSM takes the form shown in Fig.~\ref{keyGraphOld}.
Let $\mu$ be a stable matching in this problem, $\mu(x)\not\in\{c,d\}$. Then $R_\mu(x)<r'_x$.
\end{lemma}

\noindent
\textbf{Proof of Lemma~\ref{key13}:} 
Let us prove this lemma ab contrario. In this case, $R_\mu(x)>r'_x+1$.

Let us first make sure that $\mu(c)=a$. Assume the contrary. Then $R_\mu(c)>1$.
Consider the triple $(c,a,x)$. It is not blocking, only if $\mu(a)=x$. 
But then $R_\mu(b)>1$. Moreover, according to assumptions of the lemma, $\mu^{-1}(a)\not\in\{c,d\}$, therefore $R_\mu(d)>1$ and $R_\mu(c)>1$.
Consequently, the triple $(c,b,x)$ is not blocking, only if $\mu(c)=b$ (see Fig.~\ref{keyGraphOld}).

Since $R_\mu(d)>2$, the triple $(d,b,x)$ is blocking.

Thus, we have proved that $\mu(c)=a$, consequently, $R_\mu(d)>1$. Moreover, by assumption, $\mu(a)\neq x$, whence we conclude that $R_\mu(a)>1$. Then the triple $(d,a,x)$ is blocking.
\qquad  $\square$

The assertion of Lemma~\ref{key1} evidently follows from proved lemmas~\ref{key11},~\ref{key111}, and~\ref{key13}.

\begin{remark}
\label{r1}
In all figures, vertices that characterize genders are colored so as to make graph edges be directed only from red vertices to green ones, from green vertices to blue ones, and from the latter to red ones. However, it is evident that one can ``shift these colors modulo~3''. For example, in the graph shown in Fig.~\ref{keyGraph}, we can recolor all blue vertices to green. Then, certainly, we should recolor all vertices that were originally green to red and color vertices~$c$ and~$d$ to blue; this does not affect the statement of Lemma~\ref{key1}.
\end{remark} 

Let the symbol~$H$ stand for the graph of 3DSMI. For any vertex $v\in V(H)$ we put
$
\rho_H(v)=\max\limits_{(v,v')\in E(H)} r(v,v').
$ 

\begin{theorem}
\label{unSolvabable}
Let $H$ be the graph of 3DSMI of size~$n$ with no stable matching. Let us use it for constructing the graph~$G$ of 3DSM in the following way. The graph $H$ is a subgraph of the graph~$G$ (with the same ranks of edges). To each vertex $x$ of the graph~$H$ we ``attach'' the corresponding copy of the graph $H_x$ shown in Fig.~\ref{keyGraph} with vertices $a_x, b_x, c_x, d_x, e_x, s_x, t_x\not\in V(H)$ (all subgraphs $H_x$ are pairwise disjoint). Moreover, let the value $r'_x$ in formulas~(\ref{ranksRestrict}) equal $\rho_H(x)+1$. Let us define ranks of the rest edges of the graph~$G$ of 3DSM arbitrarily. Then 3DSMI with the graph $G$ has no stable matching. Here the size of 3DSM equals~$8n$.
\end{theorem}

\noindent
\textbf{Proof}: 
Assume the contrary, i.e., assume that for 3DSM with the graph~$G$ there exists a stable matching $\mu_G$.  
We intend to construct the matching $\mu_H$ for 3DSMI defined by the graph~$H$ from the matching $\mu_G$. To this end, we will make use of Lemma~\ref{key1}.
Since the matching $\mu_H$ is not stable, we can find for it a blocking triple $(v,v',v'')$ composed of vertices of the subgraph~$H$.
Let us prove that the same triple $(v,v',v'')$ is blocking for $\mu_G$. 

Let us perform the proof adhering to the above scheme. By assumption, for any vertex $x\in V(H)$ the assertion of Lemma~\ref{key1} is valid. In other words, either alternative~A takes place, while $y=\mu_G(x)\not\in V(H)$ and $\mu_G(y)\not\in V(H)$, or alternative~B takes place, while $(x,y)\in E(H)$ and $(y,\mu_G(y))\in E(H)$ (the violation of the latter correlation, i.e., the fulfillment of alternative~A for the vertex~$y$, leads to a contradiction, namely, $\mu_G(\mu_G(y))\in\{a_y,b_y,t_y\}$. Therefore $\mu_G(\mu_G(y))\not\in V(H)$, but by the initial assumption, $\mu_G(\mu_G(y)) =x$ is a vertex of the subgraph $H$).

Note that if alternative~A takes place, then
\begin{equation}
\label{eneq}
\rho_H(x)<R_{\mu_G}(x).
\end{equation}

Let us associate the matching $\mu_G$ with the matching $\mu_H$ of 3DSMI with the graph~$H$. Assume that in the case of alternative~A, $\mu_H(x)=x$ (i.e., the agent $x$ remains single). In the case of alternative~B, we put $\mu_H(x)=\mu_G(x)$.

By condition, the matching $\mu_H$ is not stable. Therefore, there exists a blocking triple $(v,v',v'')$, where $v,v',v''\in V(H)$. The blocking property means that each vertex in this triple satisfies the inequality that connects ranks, for example, for the vertex~$v$ this inequality takes the form $r(v,v')<R_{\mu_H}(v)$. If $\mu_H(v)\neq v$, then $\mu_G(v)=\mu_H(v)$ and, consequently, $r(v,v')<R_{\mu_G}(v)$ 
(because in the case under consideration, $R_{\mu_G}(v)=R_{\mu_H}(v)$). 

But if $\mu_H(v)=v$, then $R_{\mu_H}(v)=\infty$. However, now we can make use of inequality~\eqref{eneq}, which implies that $r(v,v')<R_{\mu_G}(v)$.

Analogous alternatives take place for vertices $v'$ and~$v''$, namely, $r(v',v'')<R_{\mu_G}(v')$ and $r(v'',v)<R_{\mu_G}(v'')$. 
Therefore, the same triple $(v,v',v'')$ is blocking for the matching~$\mu_G$ in 3DSM.

Let us calculate the size of 3DSM. Each triple of vertices of the graph~$H$, which are associated with different genders, corresponds to three various subgraphs in the form shown in Fig.~\ref{keyGraph}; see Remark~\ref{r1} for the principle of their coloring. Note that due to the cyclic shift of colors, any such a triple of vertices in the graph $H$ is supplemented with seven new vertices of each gender of the graph~$G$; 
in other words, the number of vertices of the graph $G$ becomes 8 times greater.                                                                                                     
\qquad  $\square$
 
\section{The Key Lemma for the further reduction of the counterexample size}
\begin{figure}[h]
\begin{center}
        \begin{tikzpicture}[->,>=stealth',shorten >=1pt,auto,node distance=3cm,thick,main node/.style={rectangle,fill=blue!20,draw,font=\sffamily\Large\bfseries}]
        \node[cblue] (f) at ( 139:2) {$f$};
        \node[cgreen] (x) at ( 319:2.7) {$x$};       
        \node[cred] (c) at (0:0) {$\mathbf c$};
        \node[cred] (d) at (0:4) {$\mathbf d$};
        \node[cgreen] (e) at (90:1) {$e$};
        \node[cblue] (a) at ( 330:1.6) {$a$};
        \node[cblue] (t) at ( 29:3.3) {$t$};
        \node[cgreen] (z) at ( 292.5:1.84776) {$z$};
        \node[cblue] (b) at ( 225:1.2) {$b$};
        \node[cgreen] (s) at (15:3) {$s$};

        \path[every node/.style={font=\sffamily\small}]
        (e) edge  []  node [] {} (c)
        (c) edge  []  node [] {} (a)
        (d) edge  []  node [] {} (a)
        (a) edge  []  node [] {} (x)
        (b) edge  []  node [] {} (x)
        (t) edge  [dashed]  node [] {} (s)
        (t) edge  [dashed]  node [] {} (s)
        (t) edge  []  node [] {} (e)
        (s) edge  []  node [] {} (d)
        (e) edge  [dashed]  node [] {} (d)
        (c) edge  [dashed]  node [right] {} (b)
        (d) edge  [dashed]  node [right] {} (b)
        (x) edge  [dotted]  node [left] {} (c)
        (x) edge  [dotted]  node [right] {} (d)    
        (b) edge  [dotted]  node [right] {4} (s)    
        (d) edge  [dotted]  node [right] {3} (t)    
        (a) edge  [dotted]  node [ left] {3} (e) 
        (b) edge  [dotted]  node [left] {3} (e)
        (z) edge  [dotted]  node [left] {} (c)
        (z) edge  [dotted]  node [right] {} (d)    
        (f) edge  []  node [] {} (e)
        (a) edge  [dashed]  node [right] {} (z)
        (b) edge  [dashed]  node [right] {} (z)
        (c) edge  [dotted]  node [left] {3} (f)
        ;        
        \end{tikzpicture}
\caption{The part of the preference graph considered in Lemma~\ref{key2}. Ranks of edges $(x,c)$ and $(x,d)$ equal $r'_x$ and $r'_x+1$, those of edges $(z,c)$ and $(z,d)$ equal $r'_z$ and $r'_z+1$, correspondingly, the rank of edges $(c,f)$, $(b,e), (a,e)$, and $(d,t)$ equals $3$, while $r(b,s)=4$.}
\label{keyGraph2}
\end{center}
\end{figure}
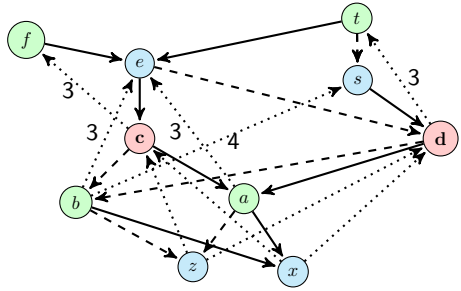

As was mentioned above, Theorem~\ref{unSolvabable}, along with the result obtained in the paper~\cite{old}, allows one to construct instances of 3DSM of size~24 with no stable matching. In the next section, we reduce the size of such counterexample to~20. As a base we use Lemma~\ref{key2}, which is a certain modification of Lemma~\ref{key1}. Recall that in Lemma~\ref{key1} we consider the subgraph (see Fig.~\ref{keyGraph}), whose copy is ``attached'' to each vertex of the graph~$H$ that defines 3DSMI with no stable matching. In Lemma~\ref{key2}, we consider the subgraph, which in certain cases can have ``two attachments'' (vertices $x$ and $z$) to two vertices of such a graph~$H$. The ``cost'' of this effect is the supplement of the subgraph with the vertex~$f$ (apart from the ``attached'' vertex~$z$).
In addition, in Lemma~\ref{key2} we impose additional constraints on preimages of vertices $x$ and $z$, which ensure the fulfillment of alternatives mentioned in this lemma. 
See Fig.~\ref{keyGraph2} for the graph under consideration.

\begin{remark}
The graph shown in Fig.~\ref{keyGraph} is a part of the graph shown in Fig.~\ref{keyGraph2}, only ranks of three edges in it are different. 
Namely, now the rank of edges directed from vertices $a$ and $b$ to the vertex~$z$ equals~2. Correspondingly, ranks of all edges that go from these vertices, which originally were not less than~2, now are larger by one, i.e., in the new graph, $r(a,e)=r(b,e)=3$ and $r(b,s)=4$. 
Ranks of all the rest edges in the subgraph of the graph shown in Fig.~\ref{keyGraph2}, which contains the same vertices~$x, a, b, c, d, e, s$, and $t$, are the same as in Fig.~\ref{keyGraph} (and no new edges appear in this subgraph). 
\label{r2}
\end{remark} 

\begin{lemma}[The Key Lemma for Theorem~\ref{exampleTh}]
\label{key2}
Assume that some subgraph~$H$ of the graph of 3DSM takes the form shown in Fig.~\ref{keyGraph2}, where ranks are indicated near the corresponding edges. Assume that $\mu$ is a stable matching such that $\mu^{-1}(z)\not\in V(H)$, provided that $R_\mu(z)<r'_z$.
Then\\ 
A) either $\mu(x)\in\{c,d\}$ and $\mu^{-1}(x)\in\{a,b,f,t\}$; \\
B) or $R_\mu(x)<r'_x$. Moreover, if, in addition, $\mu^{-1}(x)\not\in V(H)$, then
\\ \phantom{ttt} a) either $\mu(z)\in\{c,d\}$ and $\mu^{-1}(z)\in\{a,b,t\}$, 
\\ \phantom{ttt} b) or $R_\mu(z)<r'_z$.
\end{lemma}

Let us prove Lemma~\ref{key2} with the help of Lemma~\ref{key13}. But first let us consecutively prove analogs of lemmas~\ref{key11} and~\ref{key111}.
\begin{figure}[h]
\begin{center}
        \begin{tikzpicture}[->,>=stealth',shorten >=1pt,auto,node distance=3cm,thick,main node/.style={rectangle,fill=blue!20,draw,font=\sffamily\Large\bfseries}]
        \node[cblue] (f) at ( 139:2) {$f$};
        \node[cred] (c) at (0:0) {$\mathbf c$};
        \node[cgreen] (e) at (90:1) {$e$};
        \node[cblue] (a) at ( 330:1.6) {$a$};
        \node[cblue] (b) at ( 225:1.2) {$b$};
        
        \node[cgreen] (x) at ( 305:2.0) {$x$};
        \node[cgreen] (z) at ( 270:1.70711) {$z$};

        \path[every node/.style={font=\sffamily\small}]
        (e) edge  []  node [] {} (c)
        (c) edge  []  node [] {} (a)
        (a) edge  []  node [] {} (x)
        (b) edge  []  node [] {} (x)
        (c) edge  [dashed]  node [right] {} (b)
        (x) edge  [line width=0.6mm,dotted]  node [left] {} (c)
        (a) edge  [dotted]  node [right] {3} (e)
        (b) edge  [dotted]  node [left] {3} (e)
        (f) edge  []  node [] {} (e)
        (a) edge  [dashed]  node [right] {} (z)
        (b) edge  [dashed]  node [right] {} (z)
        (c) edge  [dotted]  node [left] {3} (f)
        ;        
        \end{tikzpicture}
\caption{The part of the preference graph considered in Lemma~\ref{key21}.}
\label{keyGraph21}
\end{center}
\end{figure}
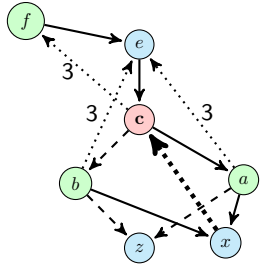

\begin{lemma}
\label{key21}
Assume that some subgraph of the graph~$G$ of 3DSM takes the form shown in Fig.~\ref{keyGraph21}, where ranks are indicated near the corresponding edges. Let $\mu$ be a stable matching such that $\mu(x)=c$. 
Then $\mu(c)\in\{a,b,f\}$.
\end{lemma}

\noindent
\textbf{Proof of Lemma~\ref{key21}:}
By assumption, $R_\mu(e)>1$. Assume that $\mu(c)\not\in \{a,b,f\}$. Then $R_\mu(c)>3$, $R_\mu(a)>1$, and $R_\mu(b)>1$. The triple $(a,e,c)$ (the triple $(b,e,c)$) is not blocking, only if either $\mu(a)=z$, or $\mu(a)=e$ (either $\mu(b)=z$, or $\mu(b)=e$). Since $\{\mu^{-1}(e),\mu^{-1}(z)\}=\{a,b\}$, we conclude that $\mu(f)\neq e$ and $R_\mu(f)>1$. But then the triple $(f,e,c)$ is blocking.
\qquad  $\square$

\begin{lemma}
\label{key211}
Assume that some subgraph~$H$ of the graph of 3DSM takes the form shown in Fig.~\ref{keyGraph2}.
Assume that $\mu$ is a stable matching in this problem, $\mu(x)=d$, and, moreover, if $R_\mu(z)<r'_z$, then $\mu^{-1}(z)\not\in V(H)$. Then $\mu(d)\in\{a,b,t\}$.
\end{lemma}

\begin{remark}
\label{r3}
In the case, when $R_\mu(z)<r'_z$, the assertion of Lemma~\ref{key211} actually coincides with that of Lemma~\ref{key111}. Really, in this case, by assumptions of Lemma~\ref{key211}, $\mu(a)\neq z$ and $\mu(b)\neq z$. Therefore, the shift of ranks of edges mentioned in Remark~\ref{r2} does not matter, because the order of these ranks remains the same and 
$\mu(z)\not\in V(H)$, $\mu^{-1}(z)\not\in V(H)$. 
\end{remark}

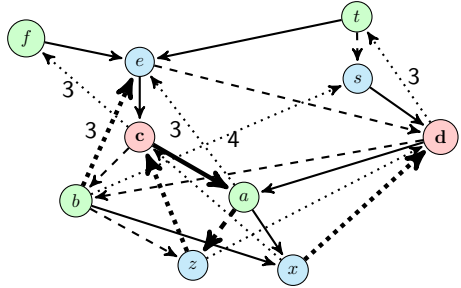
\begin{figure}[h]
\begin{center}
        \begin{tikzpicture}[->,>=stealth',shorten >=1pt,auto,node distance=3cm,thick,main node/.style={rectangle,fill=blue!20,draw,font=\sffamily\Large\bfseries}]
        \node[cblue] (f) at ( 139:2) {$f$};
        \node[cgreen] (x) at ( 319:2.7) {$x$};       
        \node[cred] (c) at (0:0) {$\mathbf c$};
        \node[cred] (d) at (0:4) {$\mathbf d$};
        \node[cgreen] (e) at (90:1) {$e$};
        \node[cblue] (a) at ( 330:1.6) {$a$};
        \node[cblue] (t) at ( 29:3.3) {$t$};
        \node[cgreen] (z) at ( 292.5:1.84776) {$z$};
        \node[cblue] (b) at ( 225:1.2) {$b$};
        \node[cgreen] (s) at (15:3) {$s$};

        \path[every node/.style={font=\sffamily\small}]
        (e) edge  []  node [] {} (c)
        (c) edge  [line width=0.6mm]  node [] {} (a)
        (d) edge  []  node [] {} (a)
        (a) edge  []  node [] {} (x)
        (b) edge  []  node [] {} (x)
        (t) edge  [dashed]  node [] {} (s)
        (t) edge  []  node [] {} (e)
        (s) edge  []  node [] {} (d)
        (e) edge  [dashed]  node [] {} (d)
        (c) edge  [dashed]  node [right] {} (b)
        (d) edge  [dashed]  node [right] {} (b)
        (x) edge  [dotted]  node [left] {} (c)
        (x) edge  [line width=0.6mm,dotted]  node [right] {} (d)    
        (b) edge  [dotted]  node [right] {4} (s)    
        (d) edge  [dotted]  node [right] {3} (t)    
        (a) edge  [dotted]  node [ left] {3} (e) 
        (b) edge  [line width=0.6mm,dotted]  node [left] {3} (e)
        (z) edge  [line width=0.6mm,dotted]  node [left] {} (c)
        (z) edge  [dotted]  node [right] {} (d)    
        (f) edge  []  node [] {} (e)
        (a) edge  [line width=0.6mm,dashed]  node [right] {} (z)
        (b) edge  [dashed]  node [right] {} (z)
        (c) edge  [dotted]  node [left] {3} (f)
        ;        
        \end{tikzpicture}
\caption{The part of the preference graph considered in Lemma~\ref{key211}. Bold lines represent edges that correspond to pairs that enter in families of the matching $\mu$ in the case, when $R_\mu(z)\geqslant r'_z$ and $\mu(d)\not\in\{a,b,t\}$.}
\label{graph3}
\end{center}
\end{figure}

\noindent
\textbf{Proof of Lemma~\ref{key211}:}
Let us restrict ourselves to considering the case, when $R_\mu(z)\geqslant r'_z$. 
By assumption, $R_\mu(x)=r'_x+1$ and $R_\mu(s)>1$.
Analogously to the proof of Lemma~\ref{key111}, we assume that $\mu(d)\not\in\{a,b,t\}$. 

We get inequalities $R_\mu(d)>3$, $R_\mu(b)>1$, and $R_\mu(a)>1$. 
The triple $(c,a,x)$ is not blocking, only if $\mu(c)=a$. 

Let us prove that the matching $\mu$ contains the family $(c,a,z)$. Really, otherwise $R_\mu(z)>r'_z+1$ and $R_\mu(a)>2$, and then the triple $(d,a,z)$ is blocking.

As a corollary, we get inequalities $R_\mu(e)>2$ and $R_\mu(b)>2$. 
Consequently, $\mu(b)=e$, otherwise the triple $(d,b,e)$ is blocking (see Fig.~\ref{graph3}). 
Hence we get the inequality $R_\mu(t)>1$. Recall that by assumption, $R_\mu(d)>3$. But then the triple $(d,t,e)$ is blocking.
\qquad  $\square$

\vspace{2mm}
\noindent
\textbf{Proof of Lemma~\ref{key2}:}
Lemmas~\ref{key21},~\ref{key211}, and~\ref{key13} imply that
$$
\text{either alternative~A takes place, or} R_\mu(x)<r'_x. 
$$
It remains to prove that in the case, when $R_\mu(x)<r'_x$ and $\mu^{-1}(x)\not\in V(H)$, 
\begin{equation}
\label{alt2}
\text{there takes place one of alternatives $a$ or $b$ of Lemma~\ref{key2}.}
\end{equation}

Let us apply the same technique as that used in Remark~\ref{r3}. In this case, when considering various variants of composition of the matching $\mu$, one should neglect families, which contain the vertex $x$ together with other vertices shown in Fig.~\ref{keyGraph2}. 
Therefore we can delete the vertex $x$ from this graph together with edges that are incident to it. Moreover, we can change ranks of all edges that go from vertices, which previously were origin points for edges directed to $x$, preserving the order of ranks. In other words, we can decrease ranks of all edges that originate from vertices $a$ and $b$ by one. This does not affect the validity of correlation~\eqref{alt2}.

If we now change the notation $z$ to $x$, then in accordance with Remark~\ref{r2}, Fig.~\ref{keyGraph2} will turn into the subgraph shown in Fig.~\ref{keyGraph}. Lemma~\ref{key1} implies that one of alternatives~a or~b takes place.
\qquad  $\square$

\section{An example of unsolvable 3DSM of size~20}
Let us make use of the concrete example of 3DSMI of size~3 with no stable matching, which is given in~\cite[Theorem~2]{old}. 
See Fig.~\ref{example} for the graph~$H$ of the considered problem. See the Appendix for the proof of the fact that 3DSMI that corresponds to this graph has no stable matching.

\begin{figure}[h]
\begin{center}
\begin{tikzpicture}[->,>=stealth',shorten >=1pt,auto,node distance=3cm,thick,main node/.style={rectangle,fill=blue!20,draw,font=\sffamily\Large\bfseries}]
        \node[cred] (5) at ( 0:1) {\bf 5};
        \node[cblue] (0) at ( 60:1) {0};
        \node[cgreen] (1) at ( 120:1) {1};
        \node[cred] (2) at ( 180:1) {\bf 2};
        \node[cblue] (3) at ( 240:1) {3};
        \node[cgreen] (4) at ( 300:1) {4};
        \node[cblue] (6) at ( 340:2.5) {6};
        \node[cgreen] (7) at ( 60:2) {7};
        \node[cred] (8) at ( 10:2.5) {\bf 8};

        \path[every node/.style={font=\sffamily\small}]
        (0) edge  []  node [] {} (1)
        (1) edge  []  node [] {} (2)
        (2) edge  []  node [] {} (3)
        (3) edge  []  node [] {} (4)
        (4) edge  []  node [] {} (5)
        (5) edge  []  node [] {} (0)

        (6) edge  []  node [] {} (4)
        (7) edge  []  node [] {} (8)
        (8) edge  []  node [] {} (0)
        
        (4) edge  [dashed]  node [] {} (8)
        (8) edge  [dashed]  node [] {} (6)
        (0) edge  [dashed] node [] {} (7)
        (4) edge  [dotted]  node [] {} (2)    
        (1) edge  [dashed]  node [] {} (5)
        (5) edge  [dashed]  node [] {} (3)
        (3) edge  [dashed] node [] {} (1)
        ;
\end{tikzpicture}
\caption{The graph~$H$ of 3DSMI of size~3 with no stable matching. 
For convenience, we numerate vertices~$H$ with numbers~$v$, $v=0,1,\ldots,8$. 
The value $v\bmod 3$ specifies the gender that corresponds to the vertex~$v$. The rank of each edge, which is represented by a solid line, equals~1. Dashed lines represent edges, whose rank equals~2. The rank of the edge~$(4,2)$ equals~3.}
\label{example}
\end{center}
\end{figure}
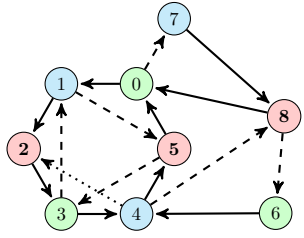

\begin{theorem}
\label{exampleTh}
Let the graph~$G$ of 3DSM contain the subgraph~$H$ shown in Fig.~\ref{example}. Assume that for all considered below subgraphs of the graph~$G$, which are copies of graphs mentioned in lemmas~\ref{key1} and~\ref{key2}, the value $r'_y$ in the corresponding copies coincides with $\rho_H(y)+1$; here $y$ is a certain vertex (we specify its number later). Thus, the graph~$G$ contains five disjoint subgraphs~$H_0$, $H_1$, $H_2$, $H_4$, and $H_7$, which are copies of the graph mentioned in Lemma~\ref{key1}; the role of the vertex $x$ is played there, correspondingly, by vertices~0, 1, 2, 4, and~7 of the graph~$H$. Subgraphs $H_0$, $H_1$, $H_2$, $H_4$, and $H_7$ have no other common vertices with the graph~$H$. Moreover, the graph $G$ has two disjoint subgraphs $H_{3,6}$ and $H_{5,8}$; they are copies of subgraphs mentioned in Lemma~\ref{key2}, the role of the vertex $x$ is played there by vertices~3 and~5, while the role of the vertex $z$ is played by vertices~6 and~8, correspondingly. 
Subgraphs $H_{3,6}$ and $H_{5,8}$ have no more common points with the graph~$H$. 
Assume that the graph~$G$ has no vertices except those considered above, i.e.,
$$V(G)= V(H_{3,6})\cup  V(H_{5,8})\bigcup_{v\in\{0,1,2,4,7\}} V(H_v).$$ We treat ranks of edges of the graph~$G$, which were not considered above, as arbitrary values. 
Then 3DSM defined by the graph~$G$ has no stable matching.
\end{theorem}

\begin{lemma}
Let~$\mu$ be a stable matching in 3DSM mentioned in assumptions of Theorem~\ref{exampleTh}. Then subgraphs $H_{3,6}$ and $H_{5,8}$ satisfy either alternative~a, or alternative~b of Lemma~\ref{key2}.
\label{last}
\end{lemma}

\noindent
\textbf{Proof of Lemma~\ref{last}:}
Let us first note that the edge $(6,4)$ is the only edge of the graph~$H$ in Fig.~\ref{example}, which originates from vertex~6. Consequently, if $R_\mu(6)<r'_6$, then $\mu(6)=4$. 
Since the subgraph $H_4$ satisfies conditions of Lemma~\ref{key1}, in this case, $\mu^{-1}(6)=\mu(4)\not\in V(H_{3,6})$, i.e., the subgraph $H_{3,6}$ satisfies conditions of Lemma~\ref{key2}.

Assume that alternative~A of this lemma takes place for the subgraph $H_{3,6}$. Then $R_\mu(3)>2$ and $R_\mu(2)>1$. The triple $(1,2,3)$ is not blocking, only if $\mu(1)=2$. But then neither alternative~A, nor~B of Lemma~\ref{key1} is possible for vertex~2. Therefore, vertex~3 satisfies alternative~B of Lemma~\ref{key2}. Then $\mu(3)\in\{1,4\}$. Since subgraphs $H_1$ and $H_4$ satisfy conditions of Lemma~\ref{key1}, we conclude that $\{\mu(1),\mu(4)\}\cap V(H_{3,6})=\emptyset$, and then $\mu^{-1}(3)\not\in V(H_{3,6})$.
Therefore, the subgraph $H_{3,6}$ satisfies conditions, which ensure the fulfillment of either alternative~a or~b of Lemma~\ref{key2}.

Assume that the subgraph $H_1$ satisfies alternative~A of Lemma~\ref{key1}, i.e., $\mu(1)\in \{c_1,d_1\}$, $\mu^{-1}(1)\in \{a_1,b_1,t_1\}$. 
Then $R_\mu(1)>2$ and $R_\mu(0)>1$. The triple $(0,1,5)$ is not blocking, only if $\mu(5)=0$. Then for vertex~0 alternative~A mentioned in Lemma~\ref{key1} is impossible and, consequently, $\mu(0)=7$, i.e., the matching $\mu$ contains the family $(7,5,0)$. Thus, both alternative~A and~B mentioned in Lemma~\ref{key1} are impossible for vertex~7. Therefore, vertex~1 satisfies alternative~B mentioned in Lemma~\ref{key1}. 

\begin{figure}[h]
\begin{center}
\begin{tikzpicture}[->,>=stealth',shorten >=1pt,auto,node distance=3cm,thick,main node/.style={rectangle,fill=blue!20,draw,font=\sffamily\Large\bfseries}]
        \node[cred] (5) at ( 0:1) {$\overline{\mathbf 5}$};
        \node[cblue] (0) at ( 60:1) {0};
        \node[cgreen] (1) at ( 120:1) {\underline 1};
        \node[cred] (2) at ( 180:1) {\bf 2};
        \node[cblue] (3) at ( 240:1) {\underline 3};
        \node[cgreen] (4) at ( 300:1) {4};
        \node[cblue] (6) at ( 340:2.5) {6};
        \node[cgreen] (7) at ( 60:2) {7};
        \node[cred] (8) at ( 10:2.5) {\bf 8};

        \path[every node/.style={font=\sffamily\small}]
        (0) edge  []  node [] {} (1)
        (1) edge  []  node [] {} (2)
        (2) edge  []  node [] {} (3)
        (3) edge  [line width=0.6mm]  node [] {} (4)
        (4) edge  []  node [] {} (5)
        (5) edge  []  node [] {} (0)

        (6) edge  []  node [] {} (4)
        (7) edge  []  node [] {} (8)
        (8) edge  []  node [] {} (0)
        
        (4) edge  [dashed]  node [] {} (8)
        (8) edge  [dashed]  node [] {} (6)
        (0) edge  [dashed] node [] {} (7)
        (4) edge  [dotted]  node [] {} (2)    
        (1) edge  [dashed]  node [] {} (5)
        (5) edge  [dashed]  node [] {} (3)
        (3) edge  [dashed] node [] {} (1)
        ;
\end{tikzpicture}
\caption{The subgraph~$H$ considered in the proof of Lemma~\ref{last} in the case, when the graph $H_{5,8}$ satisfies alternative~A.
We underline vertices $x$ such that $(x,\mu(x))\in E(H)$ and overline the rest ones. The solid line that represents the edge $(3,4)$ illustrates the fact that $\mu(3)=4$.
}
\label{example2}
\end{center}
\end{figure}
Note that if $R_\mu(8)<r'_8$, then $\mu(8)\in\{6,0\}$. Since the subgraph $H_{3,6}$ (the subgraph $H_0$) satisfies alternative~a or~b of Lemma~\ref{key2} (alternative~A or~B of Lemma~\ref{key1}), in this case, $\mu^{-1}(8)\not\in V(H_{5,8})$. Therefore, the subgraph $H_{5,8}$ satisfies conditions of Lemma~\ref{key2}.

Assume that the subgraph $H_{5,8}$ satisfies alternative~A of this lemma. Then $R_\mu(5)>2$ and $R_\mu(4)>1$. The triple $(3,4,5)$ is not blocking, only if $\mu(3)=4$ (see Fig.~\ref{example2}). But the equality $\mu(1)=5$ is impossible for vertex~1, because $\mu^{-1}(5)\in \{a_5,b_5,f_5,t_5\}$. Therefore, $\mu(1)=2$, and vertex~2 satisfies alternative~B of Lemma~\ref{key1}. Then the matching $\mu$ contains the family $(1,2,3)$, which contradicts the proposition that $\mu(3)=4$. Consequently, the subgraph $H_{5,8}$ satisfies alternative~B of Lemma~\ref{key2}. Since $\mu(5)\in\{0,3\}$, we conclude that $\mu^{-1}(5)\not\in V(H_{5,8})$. Therefore, the subgraph $H_{5,8}$ satisfies conditions, which ensure the fulfillment of alternative~a or~b of Lemma~\ref{key2}.
\qquad  $\square$

\vspace{2mm}

\noindent
\textbf{Proof of Theorem~\ref{exampleTh}:} 
Assume that in considered 3DSM there exists a stable matching $\mu_G$. According to Lemma~\ref{last}, each vertex $x\in V(H)$ satisfies one of the following two alternatives: either $\mu_G(x)\not\in V(H)$ and $\mu_G^{-1}(x)\not\in V(H)$, or $(x,\mu_G(x))\in E(H)$ (in particular, the latter case takes place for vertices $x=1,3,5$). 
Evidently, in the case of the second alternative, $(y,\mu_G(y))\in E(H)$ for $y=\mu_G(x)$ (see the same correlation in the proof of Theorem~\ref{unSolvabable}).

Let us associate the stable matching $\mu_G$ with the matching $\mu_H$ of 3DSMI with the graph~$H$. Let us do it similarly to the proof of Theorem~\ref{unSolvabable}, namely, in the case, when the first alternative takes place for $x\in V(H)$, we put $\mu_H(x)=x$, while in the case, when the second alternative takes place, we do $\mu_H(x)=\mu_G(x)$.

There exists no stable matching for the graph~$H$, therefore, for the matching $\mu_H$ one can find a blocking triple $(v,v',v'')$, where $v,v',v''\in V(H)$. Similarly to the proof of Theorem~\ref{unSolvabable}, we conclude that the same triple $(v,v',v'')$ is blocking for the matching~$\mu_G$ of 3DSM. 
\qquad  $\square$

Let us calculate the size of 3DSM mentioned in Theorem~\ref{exampleTh}. 
Subgraphs $H_0$, $H_2$, and $H_7$ have, in total, 8 vertices, which correspond to various genders. 
Each of subgraphs $H_{3,6}$, $H_{5,8}$, together with 
the subgraph shown in Fig.~\ref{keyGraph2} (with the same coloring of vertices as that shown in the figure), has 10 vertices of each of three genders.
However, instead of the subgraph shown in Fig.~\ref{keyGraph2}, we get subgraphs $H_1$ and~$H_4$, i.e., two identical subgraphs shown in Fig.~\ref{keyGraph} (with the same coloring of vertices). They have, in total, two vertices of each gender, which are ``extra'' in comparison with vertices of the subgraph shown in Fig.~\ref{keyGraph2}. 
Therefore, the graph~$G$ mentioned in Theorem~\ref{exampleTh} defines 3DSM of size~8+10+2=20.

\section{Concluding remarks and open problems}
In this paper, we decrease the size of an instance of 3DSM with no stable matching. The size of the initial example proposed by C.K.~Lam and C.G.~Plaxton equals $n=90$; the size of the example of the same problem given by us here equals $n=20$.
Earlier K.~Pashkovich and L.~Poirrier (\cite{new}) have proved that in any 3DSM of size $n\leqslant 5$ there exists a stable matching. Therefore, the minimum size of 3DSM with no stable matching is not less than 6 and not greater than 20. Its exact value is not known yet. Possibly, computer experiments with unsolvable 3DSMI, including that of the minimum size stated by us earlier in~\cite{old}, can be helpful in answering this question.

B.~Pittel has pointed out that the average value of stable matchings grows faster than $n^2\ln^2(n)$ as $n\to\infty$ (\cite{Pittel2}). 
Therefore, it seems natural to suppose that the percentage of unsolvable 3DSM tends to zero as $n\to\infty$. The probabilistic technique~\cite{alon}, which is often used for counting the number of solutions to combinatorial problems, seems to be useful in studying the mentioned problem.

Naturally, it is interesting to obtain a generalization of the classical 2-dimensional stable matching problem (2DSM) to the case of three genders without such exceptions. 
As is well known, stable matchings always exist in 3GSM with V.I.~Danilov lexicographically acyclic preferences~\cite{danilov}. In this case, genders are not symmetric; a hierarchy of their ``importance’’ is assumed to be given a~priori. However, we can state a symmetric variant of the problem, whose solutions also exist (and their quantity is larger than that in the non-symmetric case). 

Let us consider this case in detail. Assume that the following six preference matrices are given: $r_{MW}$ is the matrix of preferences of women among men, $r_{MD}$ is the matrix of preferences of dogs among men; the sense of analogous denotations $r_{WM}$, $r_{WD}$, $r_{DM}$, and $r_{DM}$ is evident. (Following the monograph by D.~Knuth, we use matrix denotations instead of order symbols accepted in the game theory.) 
For denoting ranks of agents~$x$ in a 3D-matching $\mu$, we need to modify the denotation $R_\mu(x)$ used by us earlier in this paper. Since now the agent $x$ interacts with two partners, we introduce one more superscript that denotes the gender of the partner.
Thus, if $(m,w,d)\in\mu$, then, for example, $R_{\mu,W}(m)=r_{MW}(m,w)$, analogously, $R_{\mu,D}(m)=r_{MD}(m,d)$, and so on. Let $(m,w,d)\not\in\mu$. We say that a triple $(m,w,d)$ is weakly blocking for a matching~$\mu$, if
\begin{eqnarray*}
&r_{MW}(m,w)\leqslant R_{\mu,W}(m),\quad &r_{WM}(w,m)\leqslant R_{\mu,M}(w),\\ 
&r_{DM}(d,m)\leqslant R_{\mu,M}(d),\quad &r_{MD}(m,d)\leqslant R_{\mu,D}(m),\\
&r_{WD}(w,d)\leqslant R_{\mu,D}(w),\quad &r_{DW}(d,w)\leqslant R_{\mu,W}(d)
\end{eqnarray*}
(see correlations~\eqref{oldc}). Note that here $m,w,d$ do not necessarily belong to three distinct triples in the matching~$\mu$, it suffices that at least one of three elements of the blocking triple ``is new’’. One can easily prove that due to the latter condition, at least two inequalities among six ones given above are strict. 
We say that a matching $\mu$ is (strong) \textit{stable}, if it contains no weakly blocking triple. 

We can easily construct such matchings, for example, in the following way. First we solve 2DSM on pairs from $M\times W$ with preference matrices $r_{MW}$ and $r_{WM}$, and do 2DSM on pairs from $M\times D$ with preference matrices $r_{MD}$ and $r_{DM}$. Then we ``combine'' solutions by transforming pairs $(m',w')$ in the first solution and pairs $(m',d')$ in the second one (i.e., pairs consisting of two solutions with the same elements from $M$) to triples $(m',w',d')$ of the matching $\mu$ of the three-dimensional problem.
One can easily see that in this case, any triple $(m,w,d)$, which does not enter in $\mu$, cannot be blocking, because otherwise one of the first four inequalities stated above is violated.

We can analogously define a new variant of 3DSMI as a problem with incomplete preference lists such that one of its solutions can be found in an evident way.

In the case of 2DSM, one can introduce a natural partial order on all its solutions (see~\cite{knuth}). J.~Conway and D.~Knuth have proved that this order forms a distributive lattice. Moreover, as was proved later~\cite{Bla,Gus}, the correspondence between distribution lattices and stable matchings is biunique. The authors of the mentioned papers have also explicitly described a partially ordered set (POSET), whose ideals define a lattice of stable matchings. Elements of this POSET represent the so-called rotations (i.e., cycles that include alternately men and women) defined by preference matrices of 2DSM~\cite{Irw} (see the monograph~\cite{GusIrw_book} for more detail about rotations and their connection with the distributive lattice). 
The structure of this lattice allows one to find a stable matching with the minimum regret for this problem within quadratic time, and to solve other optimization problems in the class of two-dimensional stable matchings within polynomial time.

We are interested in finding an analogous structure for the three-dimensional problem defined above. However, this is not so easy, because in this problem the existence of some two-dimensional projection of a solution to 3DSM, which solves the corresponding 2DSM, is not guaranteed. 

Let us give a concrete example of this problem. Let $n=6$. Assume that with $\ell=1,2$ any pair $(x,y)$ such that $r(x,y)=\ell$ satisfies the equality $r(y,x)=\ell$ (``love and strong sympathies are always mutual’’). Consider a graph with 18 vertices (all agents). Let us construct in this graph edges $(x,y)$ such that $r(x,y)=1$. Let the resulted graph represent a union of three cycles of length~6. Moreover, if we supplement this graph with edges $(x,y)$ such that $r(x,y)=2$, then all of them will lie inside the connectivity component formed by cycles. Evidently, such construction is possible and (accurate to isomorphism) uniquely defines elements of ranks~1 and~2 in matrices $r_{MW}$, $r_{MD}$, $r_{WM}$, $r_{WD}$, $r_{DM}$, and $r_{DW}$). We prove that independently of the rest ranks in these matrices, this problem has a stable matching, which consists, evidently, of 6 triples and has the following property: by deleting representatives of any fixed gender from these triples one cannot get a solution to the corresponding~2DSM.

Let us explain this property in more detail. For convenience, denote elements of three cycles mentioned above as $(m_{i,0},w_{i,0},d_{i,0},m_{i,1},w_{i,1},d_{i,1})$, $i=1,2,3$. Then one can easily see that the only solution to 2DSM with matrices $r_{MW},r_{WM}$ takes the form $(m_{i,k},w_{i,k})$, where $i=1,2,3$ and $k=0,1$. In other words, the unique solution to 2DSM consists of ``loving couples’’. Really, if some pair does not enter in the solution, then it is blocking. Analogously, the unique solution to 2DSM with matrices $r_{WD},r_{DW}$ takes the form $(w_{i,k},d_{i,k})$, $i=1,2,3$ and $k=0,1$; while the unique solution to 2DSM with matrices $r_{DM},r_{MD}$ takes the form $(d_{i,k},m_{i,(k+1)\bmod 2})$, $i=1,2,3$ and $k=0,1$. 

At the same time, 3DSM under consideration has the following solution: 
\begin{eqnarray*}
\mu &=& \{(m_{1,0},w_{1,0},d_{1,0}),(m_{1,0},w_{1,0},d_{1,0}),\\ 
& &(w_{2,0},d_{2,0},m_{2,1}),(w_{2,1},d_{2,1},m_{2,0}),\\
& &(d_{3,0},m_{3,1},w_{3,1}),(d_{3,1},w_{3,0},m_{3,0})\}.
\end{eqnarray*}
The absence of the blocking triple for it is due to the fact that for any $x$ in this solution $R_\mu(x)\leqslant 2$. Consequently, the blocking triple can belong only to one connectivity component of the graph defined above, where the matching of size~2 is stable.

Therefore, 3DSM considered here always has a solution; however, the structure of all solutions to this problem is rather complicated. We are going to proceed studying this issue, including related optimization problems, on the solution set of considered 3DSM.

\section*{Appendix. An instance of the 3DSMI of size $n=3$ with no stable matching}

\begin{lemma}[\cite{old}, Theorem~2]
The 3DSMI with the graph $G$ shown in Fig.~\ref{example} has no stable matching.
\label{ap}
\end{lemma}

\noindent
\textbf{Proof:}
There exist 7 families that form matchings in this problem, namely, $(0,1,5)$, $(0,7,8)$, $(1,2,3)$, $(1,5,3)$, $(2,3,4)$, $(3,4,5)$, and $(4,8,6)$.

Recall that a matching $\mu$ in the 3DSMI defined by the graph~$G$ is said to be \textit{complementable}, if there exists a triple of vertices $(v,v',v'')$ such that $\mu(v)=v$, $\mu(v')=v'$, $\mu(v'')=v''$, and $\{(v,v'),(v',v'')(v'',v)\}\subseteq E(G)$.

Evidently, any complementable matching is not stable, 
the triple $(v,v',v'')$ mentioned in the above paragraph is blocking for it. Therefore, for proving the absence of a stable matching, it suffices to find blocking triples for all noncomplementable matchings. 
For the graph shown in Fig.~\ref{example} there exists 8 noncomplementable matchings. Below we give their complete list together with blocking triples:\\
1) $\{(0,1,5),(2,3,4)\}$, the blocking triple is $(4,8,6)$;\\
2) $\{(0,1,5),(4,8,6)\}$, the blocking triple is $(1,2,3)$;\\
3) $\{(0,7,8),(1,2,3)\}$, the blocking triple is $(3,4,5)$;\\
4) $\{(0,7,8),(1,5,3)\}$, the blocking triple is $(2,3,4)$;\\
5) $\{(0,7,8),(2,3,4)\}$, the blocking triple is $(0,1,5)$;\\
6) $\{(0,7,8),(3,4,5)\}$, the blocking triple is $(0,1,5)$;\\
7) $\{(1,2,3),(4,8,6)\}$, the blocking triple is $(0,7,8)$ or $(3,4,5)$;\\
8) $\{(1,5,3),(4,8,6)\}$, the blocking triple is $(0,7,8)$. \qquad  $\square$
\end{document}